\documentclass[11pt]{amsart}
\usepackage{geometry}                
\usepackage[parfill]{parskip}    
\usepackage{graphicx}
\usepackage{amssymb}
\usepackage{epstopdf}
\usepackage{color}

\newtheorem{theorem}{Theorem}
\newtheorem{proposition}[theorem]{Proposition} 
 
\newtheorem{corollary}[theorem]{Corollary} 
 
\newtheorem{algorithm}[theorem]{Algorithm} 

\title{On the lexicographic representation of numbers}
\author{Vincenzo Manca\\{\tiny University of Verona - Italy}}

\begin{document}
\maketitle
{\bf Abstract}
It is proven that, contrarily to the common belief, the notion of zero is not necessary for having positional representations of numbers. Namely, for any positive integer $k$,  a positional representation with the symbols for $1, 2, \ldots , k$ is given that retains all the essential properties of the usual positional representation of base $k$ (over symbols for $0, 1, 2 \ldots , k-1$). Moreover, in this zero-free representation, a sequence of symbols identifies the number that corresponds to the order number that the sequence has in the ordering where shorter sequences precede the longer ones, and among sequences of the same length the usual lexicographic ordering of dictionaries is considered. The main properties of this lexicographic representation are proven and conversion algorithms between lexicographic and classical positional representations are given. Zero-free positional representations are relevantt in the perspective of the history of mathematics, as well as, in the perspective of emergent computation models, and of unconventional representations of genomes.

\

{\bf Keywords}\\
Positional representation, Zero, Lexicographic ordering

\section{Number Positional Representations}
The positional representation of numbers was an epochal discovery that is a landmark in the history of mathematics and science. Il is based on a notational system of Hindu origin that, via the Arabian mathematics, was imported in Europe with the famous book {\it Liber Abaci} published in 1202 by Leonardo Fibonacci (Leonardo from Pisa). This number representation was the seed of a new mathematical perspective, algebraic and algorithmic, that surely contributed to the important discoveries of Renaissance mathematics, from which modern mathematics stemmed.

An important aspect of positional representation was the notion of zero (from an arabian root common to the noun {\it zephyr}, a gentle wind whose origin is hardly determinable). 
Any natural number can be univocally represented as sum of powers of any base $k>1$ multiplied by coefficients smaller than $k$. From this unicity it easily follows that, given an alphabet of $k$ symbols (standing for $0, 1, 2, , k-1$), every number is univocally represented by a sequence of symbols in this alphabet, where the position of the symbol encodes a power of the basis $k$, and the symbol in that position encodes a multiplicative coefficient, in such a way that the sum of these multiplied powers yields the number. When a (symbol for) zero occurs, then the corresponding power gives a null contribution to the overall sum providing the number. Classical texts in number theory and history of mathematics usually connect intrinsically the notion of zero with the positional notation \cite{ore,number-book}. For example, in \cite{ore} it is written that ``The only complication which the positional notation involves lies in the necessity of introducing a {\it zero symbol} to express a void or missing class; for instance, 204 is different from 24''. The same ``necessity'' is expressed by \cite{number-book}, and certainly the attitude toward this conviction was enforced by the decimal notation and the related algorithms for basic operations (using alignment, carry, positional shift, and decimal point), after the seminal works of Fran\c{c}ois Vi{\`e}te and Simon Stevin, in sixteenth century \cite{chabert}. However, in \cite{toeplitz}, in the context of the historical analysis of number representation systems, some remarks clearly suggest that positional notations were present in some ancient number systems. One of them is the Greek one, especially in the sexagesimal notation used in Ptolemy's Almagest, where some kind  positional concepts, made possible the complex computations required in the contexts of geometrical and astronomical problems. In \cite{knuth} an illuminating synthesis of the historical development of number representation is reported.
In passing, we recall that, in a geographically far culture, and chronologically antecedent, Maya developed a positional number system, with base twenty, also having a well defined notion of zero  \cite{zero}.

In this  note we show that the connection between zero and the positional notation is not necessary. Namely, we prove that the usual lexicographic representation of strings, in a given alphabet, is a positional number system, but without  any symbol for zero.
The essence of positional notations is the notion of cipher sequence, where ciphers have positional values corresponding to the powers of the base. The main advantage of this notational mechanism is that arithmetical operations (firstly sum and multiplication) can be calculated by {\bf operation tables}, that is, by only knowing  the values of the given operation when its arguments range over a finite set of numbers (whose cardinality is equal to the value of base). We will show that just this mechanism is behind the lexicographic representation by strings, where the number $n$ is represented by the string in position $n$ according to the lexicographic enumeration. In \cite{zeroless} the possibility of a {\it zeroless}  positional system, and its connection with the string ordering, was already investigated. However, the analysis developed in \cite{zeroless} is very different from the one we present in this paper, and the main properties of a lexicographic representation are determined in a completely different  way. Moreover,  the representation algorithms are here defined in a more direct  manner (without passing through the classical positional systems as in \cite{zeroless}), and the correctness proofs are formally established (in \cite{zeroless} they are  informally motivated). Finally, here, translation functions between lexicographic and classical positional representations are formally defined. It is interesting to remark that the main idea of the present paper arose during the search for unconventional genome representation, where DNA sequences are viewed as numbers (in base four). In that context the following question was posed: {\it Define a function that assigns to any string, over a given alphabet, its order in the lexicographic enumeration}. The solution to this question is given by Equation \eqref{X1} of the ``Lexicographic Theorem" of Section 2. From the recursive formula  \eqref{X1},  the iterative formula of Equation \eqref{X2} easily follows, which directly provides a number positional representation with the properties discussed along the paper.

\section{The Lexicographic Number Representation}
Let $S = \{a_1, a_2, \ldots , a_k\}$ a finite alphabet of symbols where we define an enumeration function 
$\omega$ of its symbols  such that
$\omega(a_i)= i$, for $1 \leq i \leq k$. We denote by $S^*$  the set of strings, that is, of finite sequences over $S$. 
In this paper we denote by $|S|$ the cardinality of the set $S$, while $|\alpha|$ denotes the length of string $\alpha$, moreover string concatenation is expressed by juxtaposition, or by $*$ when it avoids confusion, and 
$\alpha(i)$ will denote the symbol of $\alpha$ in position $i$ (the first position is $1$ and the last equals the length of the string).

We can define a strict order over $S^*$ by the following conditions:
\begin{align*}
\alpha < \beta &\text{   \hspace{1cm}     if     \hspace{1cm}       } |\alpha| < |\beta|\\
\alpha < \beta  & \text{     \hspace{1cm}      if      \hspace{1cm}      } |\alpha| = |\beta| \text{    and    }  \exists j : \forall i <j : \alpha(i) = \beta(i)  \text{    and    }  \omega(\alpha(j)) < \omega(\beta(j)).
\end{align*}
The second condition is the usual criterion adopted in the enumeration of the items in a lexicon. The first one is added to the second, in order to obtain a linear ordering over (the infinite set of) all strings over the alphabet.
For example, given the (ordered) alphabet $\{A < C < G < T\}$ (typical of genomes), then Table \ref{ex1} gives the first 40 strings  over this alphabet enumerated in the lexicographic order:
\begin{table}
\begin{center}
\begin{tabular}{|l|}
\hline
A, C, G, T, AA, AC, AG, AT, CA, CC, CG, CT, GA, GC, GG, GT, TA,  \\
TC, TG, TT, AAA,  AAC, AAG, AAT, ACA, ACC, ACG, ACT, AGA, \\
AGC, AGG, AGT, ATA, ATC, ATG, ATT, CAA, CAC, CAG, CAT.\\
\hline 
\end{tabular}
\end{center}

\

\caption{The first 40 strings over $\{A,C,G,T\}$  in the lexicographic order.}
\label{ex1}
\end{table}%

\begin{theorem}[The lexicographic Theorem]
Let  $S$ be an (ordered) finite alphabet. If $x \in S$, and $\alpha \in S^*$, then  Equation \eqref{X1} inductively defines the enumeration number $\omega_S(\alpha)$ of  string 
$x\alpha$, shortly indicated by $\omega(\alpha)$, according to the lexicographic ordering on  $S^*$ ( $\omega(\lambda) =0$, if $\lambda$ is the empty string):
\begin{equation}\omega(x\alpha) = \omega(x)|S|^{|\alpha|} + \omega(\alpha) \label{X1}\end{equation}
moreover:
\begin{equation}\omega(\alpha) = \sum_{i=1}^{|\alpha|}\omega(\alpha(i))|S|^{|\alpha|-i}\label{X2}\end{equation}
\end{theorem}
{\bf Proof}. In fact, let us denote by $\omega_{=}(\alpha)$ the  position of $\alpha$ when we lexicographically enumerate the  strings of length $\alpha$ ($a^{|\alpha|}$ is the first string in position 1, if $a$ is the first symbol of $S$). In order to evaluate the position of $x\alpha$ in the lexicographic enumeration, we observe that:\\
i) $x\alpha$ follows all the strings of length $\leq |\alpha|$, which are in number:
$$\sum_{i=1}^{|\alpha|}|S|^i$$
ii) and $x\alpha$ follows the strings of length $|\alpha|+1$ beginning with any symbol $y$ such that $\omega(y) < \omega(x)$, which are in number:
$$[\omega(x)-1]|S|^{|\alpha|}$$
iii) between the first string of length $|\alpha|+1$ that begins with $x$  and $x\alpha$ (extremes included) there are 
$\omega_{=}(\alpha)$ strings.
Therefore, the position of $x\alpha$ in the lexicographic enumeration is:
$$\omega(x\alpha) = \sum_{i=1}^{|\alpha|}|S|^i + [\omega(x)-1]|S|^{|\alpha|} + \omega_{=}(\alpha)$$
that can be rewritten as:

$$ \sum_{i=1}^{|\alpha|-1}|S|^i +|S|^{|\alpha|} + [\omega(x)-1]|S|^{|\alpha|} + \omega_{=}(\alpha)$$

\noindent but, summing the two terms in the middle of the sum above, we get the term $\omega(x)|S|^{|\alpha|}$, while summing the other two (extremal) terms we get $\omega(\alpha)$. This concludes the proof of Equation \eqref{X1} (and of the inductive definition of $\omega$). The proof of Equation \eqref{X2} follows directly, by iteratively applying Equation  \eqref{X1}.
If  $k$ is the cardinality of $S$, then Equation \eqref{X2} can be rewritten as follows:  
$$\omega(\alpha) = \sum_{i=1}^{|\alpha|}\omega(\alpha(i))k^{|\alpha|-i} \text{\hspace{2cm}{\emph QED.}}$$  

An interesting simple consequence of the theorem above is the following corollary, due to the fact that, being the lexicographic order a total ordering over all the strings, every string represents the number of its lexicographic enumeration.
 
\begin{corollary} \label{unicity}
Given a number $k>0$, any non-null natural number $n$ is univocally representable by a linear combination of powers $k^i$, $c_0k^0+c_1k^1 + \ldots + c_jk^j$ where  $0< c_i \leq k$ for all $0\leq i \leq j$.
\end{corollary}

We remark that according to the corollary above, lexicographic number representation results a {\it full positional representation} because all the powers of $k$ less than a maximum value give a positive contribution to the sum representing the number, while classical representations admit the possibility of having empty contribution (expressed by zero ciphers). The lexicographic representation is the most economical way to represent the set of the first $m$ numbers. In other words, if the cost of representing these numbers is defined (with respect to a base $k$) as the sum of the lengths of strings representing them, then no representation exists that has a lower representation cost.

Table \ref{table-sum} defines operation $+$ for symbols $A,C,G,T$ that is coherent with the lexicographic ordering, as illustrated by Table \ref{gatto}. In fact, by using Table \ref{table-sum} we can  coherently add strings, in perfect accordance with their lexicographic order. 

\begin{table}
\begin{center}
\begin{tabular}{|cccccc|}
\hline
\mbox{}& + & A & C & G & T\\
\mbox{} &A & C & G & T & AA\\
\mbox{} &C & G & T & AA & AC\\
\mbox{} &G & T & AA & AC & AG\\
\mbox{} &T & AA & AC & AG & AT\\
\hline
\end{tabular}
\end{center}

\

\caption{Addition table of number lexicographic representation (four symbols).}
\label{table-sum}
\end{table}%

\begin{table}
\begin{center}
\begin{tabular}{|lcl|}
\hline
$\omega(CAT)$&=&$40$\\
$\omega(GATT)$&=&$228$\\
$CAT + GATT$ &= &$GTCT$\\
$\omega(GTCT)$ &= &$40+228=268$\\
\hline
\end{tabular}
\end{center}

\

\caption{The additivity of the lexicographic order.}
\label{gatto}
\end{table}

In general the additive table of a lexicographic system can be obtained in accordance to to the following proposition that follows directly from the previous discussion. In the following, when we want to avoid confusion, the symbols of a $k$-lexicographic representation (in base $k$) are denoted by $[1], [2], \ldots [k]$ (with numbers in decimal notation inside brackets), while symbols  of a classical $k$-positional representation (with 0) are denoted by $[0], [1], \ldots [k-1]$.

\begin{proposition}\label{prop2}
Given a lexicographic system of base $k$, 
its additive table $L_k(+)$ can be obtained from the additive table $P_{k}(+)$ of the classical positional table: i) by removing all the elements corresponding to sums $[j]+[0]$ or  $[0]+[j]$, for $j <k$; ii) by replacing in $P_{k}(+)$ any sum result $[1][0]$ by $[k]$; and iii) by adding as results of the sums $[j]+[k]$ or  $[k]+[j]$, for $j \leq k$, the elements $[1][j]$. 
\end{proposition}

The following algorithm provides the string that lexicographically represents the number $n$ over an alphabet of $k$ symbols.
It provides the inverse information with respect to the theorem above. In fact, given a positive integer $n$,  it solves the equation of unknown $\alpha$: 
$$\omega(\alpha)=n.$$

\begin{algorithm}[Lexicographic Number Representation]
For any positive integer $h$, let 
$$maxlex(k,h)= \sum_{i=1}^{h}k^{i}$$ 
$$minlex(k,h)= \sum_{i=1}^{h}k^{i-1}.$$ 
The first number expresses the greatest number lexicographically representable by strings of length  
$h$, while the second number expresses the smallest number lexicographically representable by strings of length $h$. 
Then, given a positive integer $n$, if $$minlex(k, h) \leq n \leq maxlex(k,h)$$ then it  is lexicographically representable by a string of length $h$, and we denote this value $h$ by $lex(k, n)$, or simply by $lex(n)$ when $k$ is implicitly understood. 
In this case, $n$ is univocally lexicographically represented by the string $\sigma_k(n)$ (over $\{[1], [2], \ldots, [k]\}$), according to the following recursive procedure (here $*$ denotes string concatenation and 
$\cdot$ the product):

\

$\sigma_k(n) = [n]$   if $n \leq k$\\
\indent  $\sigma_k(n) = [m] * \sigma_k(n - m\cdot k^{lex(n)-1})$ otherwise\\
\indent where $m = max \{q \leq k \; | \; n - (q \cdot k^{lex(n)-1}) \geq minlex(k, lex(n)-1)\}.$

\hspace{13.5cm}{\emph QED.}

\end{algorithm}

For example, let us consider $37$ and $36$, of course  $lex(4, 36)=lex(4, 37) = 3$, but $37$ is lexicographically represented (in the base $4$) by:
$$37 = 2\cdot 4^2 + 1\cdot 4^1 + 1\cdot 4^0$$
because $37 - 2\cdot 4^2 =5$ and $5 \geq minlex(4,2)$, while:
$$36 = 1\cdot 4^2 + 4\cdot 4^1 + 4\cdot 4^0$$
because $36 - 2\cdot 4^2 = 4$ and $4 < minlex(4,2)$.

\begin{corollary}
Given a positive integer $n$, the previous algorithm correctly provides the string $\sigma_k(n)$, representing $n$ in the lexicographic system of base $k$, such that:
$$\omega(\sigma_k(n))=n.$$
\end{corollary}
{\bf Proof}. The asserted correctness easily follows, by induction, from the definition of 
$lex$, $minlex$, $maxlex$, and the recursive procedure defining $\sigma_k(n)$. 
\hspace{4cm}{\emph QED.}

\

The correctness of a number representation $\rho$ implies that any (binary) operation $\odot$ on numbers, when performed on number representations, has to satisfy the {\it morphism condition} (in the algebraic sense), that is, $\rho$ has to commute with the operations (on numbers and representations, for the sake of simplicity both denoted by $\odot$):
\begin{equation} \rho(n) \odot \rho(m) = \rho(n\odot m).  \label{morphism} \end{equation}

Multiplication in lexicographic representation can be computed with the usual algorithm of classical positional systems, by using the multiplication Table \ref{tableX}. For example, the multiplication $37 \times 3X=147X$ can be obtained with the usual multiplication algorithm of the positional representation 
($147X$ corresponds to 1480 in the usual decimal system, coherently to the correspondence of $(37, 3X)$ to the decimal pair $(37, 40)$ and to the decimal multiplication $37 \times 40 = 1480$.

\begin{table}
\begin{center}
\begin{tabular}{|ccccccccccc|}
\hline
$\times$ & 1 & 2 & 3 & 4&5&6&7&8&9&X\\
1 &1 &2 & 3 & 4 & 5& 6& 7& 8& 9& X\\
2 &2 &4 & 6 & 8 & X& 12& 14& 16& 18& 1X\\
3 &3 &6 & 9 & 12 & 15& 18& 21& 24& 27& 2X\\
4 &4 &8 & 12 & 16 & 1X& 24& 28& 32& 36& 3X\\
5&5 &X & 15 & 1X & 25& 2X& 35& 3X& 45& 4X\\
6&6 &12 & 18 & 24& 2X& 36& 42& 48& 54& 5X\\
7&7 &14 & 21 & 28& 35& 42& 49& 56& 63& 6X\\
8&8 &16 & 24 & 32& 3X& 48& 56& 64& 72& 7X\\
9&9 &18 & 27 & 36& 45& 54& 63& 72& 81& 8X\\
X &X &1X & 2X & 3X & 4X& 5X& 6X& 7X& 8X& 9X\\
\hline
\end{tabular}
\end{center}

\

\

\caption{The multiplication table of the lexicographic representation in base ten, where the first nine symbols $1,2,3,4,5,6,7,8,9$ have the same value as in the classical decimal system, while $X$ denotes ten (zero is not present).}
\label{tableX}
\end{table}

Analogously to the additive table, the multiplicative table of a lexicographic system can be obtained in accordance to the following proposition. 

\begin{proposition} \label{prop1}
Given a lexicographic system of base $k$,
its multiplicative table $L_k(\times)$ can be obtained from the multiplicative table $P_{k}(\times)$ of the classical positional table  i) by removing all the elements corresponding to multiplications $[j]\times[0]$ or  $[0]\times[j]$, for $j <k$; ii) by replacing in $P_{k}(\times)$ any multiplication result $[1][0]$ by $[k]$, and any multiplication result $[j][0]$ by $[j-1][k]$;  iii) by adding as results of the multiplications $[j]\times[k]$ or  $[k]\times[j]$, for $j \leq k$, the elements $[j-1][k]$. 
\end{proposition}

\begin{corollary} \label{XXX}
Given a positive integer $n$, in the lexicographic system of base $k$ (symbols $[1], [2], \ldots [k]$) the result of multiplication $\sigma_k(n)\times[k]$ is represented by $\sigma_k(n-1)[k]$ ($\sigma_k(0)$ is the empty string).
\end{corollary} 

\begin{table}
\begin{align*}
&423\times \\
&8X  =\\
\_\_&\_\_\_\_\_\_\_\_\_\_\\
& 4\;2\;2\;X\\
3\;&3\;8\;4\\
\_\_\_\_&\_\_\_\_\_\_\_\_\_\_\\
3\;&7X6\;X
\end{align*}
\caption{A multiplication in the lexicographic representation of base ten. In the corresponding zero positional system  the same multiplication becomes $423 \times 90$, which usually is performed by adding a final $0$ to $423 \times 9=3807$. Of course $38070$ transforms, according to Algorithm \ref{XX}, into $37X6X$ (and {\it viceversa}). It seems that the usual multiplication is much simpler than the lexicographic one. However, also in the lexicographic system it is possible to get the result in a shorter way. In fact, the multiplication $423 \times X$ can be  obtained, according to Corollary \ref{XXX}, by appending $X$ to the predecessor of $423$, getting $422X$, then $423 \times 8=3384$, therefore, shifting and summing the two partial results, the final result is obtained.}
\label{tableXX}
\end{table}

\section{Lexicographic and ZeroPositional representations}
Rules can be established that transform a lexicographic representation, in a given base, into a classical (with zero) positional representation in the same base and {\it vice versa}. Let us denote by $\delta_k(n)$ the $k$-positional representation with zero of the number $n$ in base $k$, while $\sigma_k(n)$ is the $k$-lexicographic representation of $n$.

\begin{algorithm}[Conversion between Lexicographic and ZeroPositional representations] \label{XX}
The $k$-positional (with zero) and $k$-lexicographic representations are mutually translated, in a 1-to-1 way, by the functions 
$\theta_{0 \to k}$ and $\theta_{k \to 0}$  computed with the following algorithm.

Given a string $\alpha$ that is a positional representation with zero in base $k$, and that does not start with the symbol $[0]$ (apart the string $[0]$ that is intended to have order number zero), then factorise $\alpha$ as a concatenation of substrings (some of $\alpha_i, \gamma_i$ can be empty strings):

$$\alpha_0\gamma_1\alpha_1\gamma_2 \ldots \alpha_{h-1}\gamma_h\alpha_{h}$$

\

where $\alpha$ factors are 0-free and $\gamma$ factors have the form $[j][0]^i$ with $j>0$ and $i>0$. Now, let:

\

 $$\gamma^*_i = ([j][0]^{i-1} - [1])[k]$$ 
 
\

 where minus operation is computed in the $k$-positional representation with zero. Then:

$$\theta_{0 \to k}(\alpha) = \alpha_0\gamma^*_1\alpha_1\gamma^*_2 \ldots \alpha_{h-1}\gamma^*_h\alpha_{h}.$$

\

Analogously, given a string $\beta$ that is a lexicographic representation in base $k$, factorise $\beta$ as a concatenation of substrings (some of $\beta_i, \eta_i$ can be empty strings):

$$\beta_0\eta_1\alpha_1\eta_2 \ldots \beta_{g-1}\eta_g\beta_{g}$$

\

where $\beta$ factors are [k]-free, and $\eta$ factors have the form $[j][k]^i$ with $j<k$ and $i>0$. Now, let:

\

$$\eta^o_i = [j+1][0]^{i}$$ 
 
 \
 
where in the case that $j+1=k$ the cipher $[k]$ has to be replaced by $[1][0]$.
 Then:

$$\theta_{k \to 0}(\beta) = \beta_0\eta^o_1\beta_1\eta^o_2 \ldots \beta_{g-1}\eta^o_g\beta_{g}.$$

\end{algorithm}
\hspace{13cm}{\emph QED.}

\

The following proposition states the correctness of the translation functions $\theta_{0 \to k}$ and $\theta_{k \to 0}$ defined by the algorithm above.

\

\begin{proposition}
$$\theta_{0 \to k}(\delta_k(n)) = \sigma_k(n)$$
$$\theta_{k \to 0}(\sigma_k(n)) = \delta_k(n)$$
\end{proposition} 
{\bf Proof}.
Firstly, we observe that if we order the alphabet $\{[0], [1],\ldots , [k-1]\}$ with the natural order (where $[0]$ is the minimum symbol), then the number represented by a string $\alpha$ is the enumeration order of this string, according to the lexicographic ordering (restricted to the strings that do not start with the symbol $[0]$). Let us denote by $\omega_0(\alpha)$ the order number of $\alpha$ in this ordering. Then,
the asserted equations are equivalent to the following ones. In other words, the translations between the two representations preserve their interpretations as numbers:
$$\omega(\theta_{0 \to k}(\delta_k(n)))=n$$ 
$$\omega_0(\theta_{k \to 0}(\sigma_k(n)))= n.$$
In fact, the translation rules of the algorithm above guarantee that translations preserve the numerical value expressed by the strings, because they essentially express $[k]$ as $[1][0]$ and {\it vice versa}, but at same time, when this substitution is applied, a carry cipher $[1]$  has to be propagated (in addition or in subtraction, respectively) to the ciphers on the left of the position where substitution was performed.
\hspace{7cm}{\emph QED.}\\

For example (here base $10$ is denoted by $X$):
$$\theta_{X \to 0}(2X9X5) =  301005.$$
$$\theta_{0 \to X}(301005) =  2X9X5$$
$$\theta_{X \to 0}(2XXX9XX5) =  300010005$$
$$\theta_{0 \to X}(300010005) =  2XXX9XX5.$$

The proof of the previous proposition puts in evidence that classical positional representations can be seen as a restricted form of lexicographic representation (limited to a subset of strings over the $k$-base alphabet). 

\begin{corollary}
Given a base $k>1$, any number has, in this base, a unique positional representation with zero.
\end{corollary}
{\bf Proof}. The unicity of lexicographic representation is obvious (see Corollary \ref{unicity}), therefore, the unicity of positional representation with zero follows immediately from it, via the translation $\theta_{k \to 0}$. \hspace{11cm}{\emph QED.}\\

\


A final remark will be useful to the discussion of the conclusive section. Even when the base of a positional system is a number greater than 10 (for example, 20 or 60), multiplication algorithms based on positional representation can be very efficient. In fact, by using some strategy of distributed carry, as the lattice (or sieve) multiplication method, originally pioneered by Islamic mathematicians,
and described by Fibonacci in his {\it Liber Abaci},  and using only a subset of the multiplication table (and/or other tricks), even complex multiplications can be easily developed. Tables \ref{symbolic} and \ref{symbolic-60} give two examples, where a kind of lattice multiplication is joined to the restriction of multiplication tables (to 2 and 5 in the first case, and to 2, 5, 10 in the second case).


\begin{table}
\begin{align*}
427\times &\\
35=&\\
\_\_\_\_\_\_\_\_\_\_\_\_\_\_\_\_\_\_\_\_\_\_\_&\_\_\_\_\_\_\_\\
(5\cdot 4)(5\cdot 2)&(5\cdot 7)\\
(3\cdot 4)(3\cdot 2)(3\cdot 7)&\\
\_\_\_\_\_\_\_\_\_\_\_\_\_\_\_\_\_\_\_\_\_\_\_&\_\_\_\_\_\_\_\\
\fbox{$0$}''''\;\; \;\;\;\;\; \;\;\; \;\; \fbox{$0$}''&\;\;\;\;\;\fbox{$5$}'\\ 
(3\cdot 2)\;\;\;\;\;\; \;\;6\;\;\;\;\;\; \;\fbox{$(3\cdot 5)$}'''_{\;5}&\\
(3\cdot 2)\;\; \; \;\;\;\;\fbox{1}''\;\;\;\,\;\;\;\; (3\cdot 2)\;\;\;&\\
\;\;\fbox{2}''''\;\;\;\;\;\fbox{1}'''\;\;\;\,\;\;\;\;\fbox{3}' \;\;\;\;\;\;\;&\\
\_\_\_\_\_\_\_\_\_\_\_\_\_\_\_\_\_\_\_\_\_\_\_\_\_\_\_\_\_\_\_&\_\_\_\_\_\_\_\\
\;\;\;\;\;&\;\;5\\
\;\;6 \;\;\;\;\;6\;\;\;\;\;5\;\;&\;\;\;\\
6\;\;\;\;\;1\;\;\;\;\; 6\;\;&\;\;\;\\
\;\;\;\;2\;\;\;\;\;1 \;\;\;\;\;3\;\;&\\
\_\_\_\_\_\_\_\_\_\_\_\_\_\_\_\_\_\_\_\_\_\_\_\_\_\_&\_\_\_\_\_\_\_\\
1\;\;\;\;\;\;4\;\;\;\;9\;\;\;\;\;4\;\;\;&\;\;5\\
\end{align*}
\caption{A multiplication in decimal representation, developed by a method of {\it distributed carry}, where only multiplications by 2 and 5 are used.
In this method, when a multiplication provides a carry, then it is appended to the column on the left. Multiplications by ciphers different from 2 or 5 are reduced to sums of them ($3\cdot 7= 3\cdot  5 + 3 \cdot  2$). Boxes and apices are used, for clarity's sake, by decorating a carry and the multiplication generating it, with the same number of apices.  The algorithm remains essentially the same for lexicographic representations (using the appropriate multiplication tables).}
\label{symbolic}
\end{table}

\begin{table}
\begin{align*}
&[7][7]\times \\
&[35]=\\
&\_\_\_\_\_\_\_\_\_\_\_\_\_\_\_\_\_\_\_\_\_\_\_\_\_\_\_\\
&([35]\cdot [7])([35]\cdot [7])\\
&\_\_\_\_\_\_\_\_\_\_\_\_\_\_\_\_\_\_\_\_\_\_\_\_\_\_\_\\
&\_\_\_\_\_\_\_\_\_\_\_\_\_\_\_\_\_\_\_\_\_\_\_\_\_\_\_\\
&([10]\cdot [7])([10]\cdot [7])\\
&([10]\cdot [7])([10]\cdot [7])\\
&([10]\cdot [7])([10]\cdot [7])\\
&([5]\cdot [7])([5]\cdot [7])\\
&\_\_\_\_\_\_\_\_\_\_\_\_\_\_\_\_\_\_\_\_\_\_\_\_\_\_\_\\
&[10][10]\\\
&[10][10]\\
&[10][10]\\
&[35][35]\\
[1]&\,\;[1]&\\
[1]&\:\;[1]&\\
[1]&\;\;[1]&\\
\_\_\_\_\_\_\_\_\_&\_\_\_\_\_\_\_\_\_\_\_\_\_\_\_\_\_\_\_\_\_\\
[4]\;&[9]\;\;\;[5]\\
\end{align*}
\caption{The same multiplication of Table \ref{symbolic} realized in sexagesimal representation, with the method of {\it distributed carry}, and where only multiplication by 10 and 5 are used (if $[j]$ is the $j$-sexagesimal cipher, then $427=[7][7], 35 = [35], 14945 = [4][9][5]$). Here the sub-notation (inside brackets) is decimal. The adoption of another representation of cyphers (for instance, that one of {\it Almagest}) is only matter of ``syntactic sugar''.}
\label{symbolic-60}
\end{table}

\section{Conclusions}
The French word for computer is {\it ordinateur}. Certainly, it relates to the fundamental arithmetic ordering generated from zero by means of successor operation: $0, 1, 2, \ldots$ as basis of any numerical calculus. From successor all the operations come via iteration: sums are successor iterations, multiplications are sum iterations, exponentials are multiplication iterations, differences iterate the predecessor operation (the inverse of successor ordering), divisions iterate difference, and so on.  What we discover from the previous discussion is a further intriguing connection between calculus and ordering. In fact, string orderings are intrinsically related to number positional representations, which provide efficient methods for computing the arithmetical operations.

In \cite{kaplan} the fundamental role of positional representation in the process of understanding numbers is explained by many specific examples and analyses. For instance, the existence of irrational numbers, a masterpiece of Greek mathematics, becomes a simple corollary of  sequence  representations of numbers. In fact, it is easy to prove that the division of two positive integers, in positional representation, provides a finite sequence of ciphers or an infinite sequence of ciphers that is periodic (it is a consequence of simple combinatorial arguments). Then, any rule generating an infinite sequence of ciphers that is not periodic has to denote a number that cannot be rational. More sophisticated examples involve results about real and transcendental numbers \cite{kaplan}. It is also worthwhile to recall here that the famous Turing's paper of 1936 \cite{turing}, which was the starting point of the theory of computability, was aimed at providing a definition of {\it computable real number}, as algorithmically generated infinite sequence of ciphers.

The kind of positional representation that is inherent to the lexicographic notation could be the basis for explaining some intriguing aspects of ancient mathematics. In fact, we know that, even before Greek mathematics, in ancient civilisations, complex calculations were developed (Babylonians were able to manage, in base 60, numerical algorithms related to complex astronomical problems).  This ability  is surely impossible without an efficient number representation. 


In \cite{toeplitz} the solution is reported that was given by Fibonacci to the cubic equation $x^3+2x^2+10x=20$ in sexagesimal notation (with an error only to the sixth fractional sexagesimal cipher ($60^{-6}$). This reveals an amazing mathematical mastery of computations with sexagesimal fractions, additions, multiplication, extraction of square roots, and so on. However, sexagesimal tradition goes back to the {\it Almagest} of Ptolemy, composed about A.D. 150, within Greek culture, where tables of chords (essentially sines) are given for every $(1/2)^o$ from $0^o$ to $180^o$, with sexagesimal primes and seconds, which is equivalent to an approximation greater than three decimals fractions (for instance, $\sin 1 = 1 '  2 '' 50 '''$). Ptolemy never showed numerical computations, therefore we  do not know to what extent he used positional principles in computing, but according to our previous analysis, we know that they are independent of zero. Moreover, as it is remarked in  \cite{toeplitz}, Ptolemy uses symbol $\mbox{}^o$ in the first position as indication of angle degree, but in the following positions (for primes, seconds, and so on) symbol  $\mbox{}^o$ is used to denote the absence of any value for that sexagesimal position, that is, what now we call zero. In conclusion, also this confirms that sexagesimal notation is a positional system, within which arithmetic efficient algorithms can be available, with the full computational power of positional notation, even without any explicit notion of zero. In passing, with Ptolemy, Greek mathematics incorporated concepts coming from oriental cultures (especially Babylonian astronomy), and, at that time, they were something well known and established. Probably, the positions of the sexagesimal system are related to the notion of time periods, as suggested by astronomical observations (already Archimedes, investigating on the representation of big numbers, elaborated a method of number representation based on periods). 
  
More details on the historical aspect are beyond the aims of our investigation (see \cite{knuth,zero,ifrah} for deeper analyses), however, surely positional systems were related to the {\it abacus} where ciphers are expressed, in some way, by the number of some items in each compartments (compartments play the role of positions of symbols). If we consider a sort of ``altimetric'' variant of the abacus, we can suggest a strict relationship between positional representation with zero and lexicographic positional representation. In fact, let us encode ciphers with the height of some object, say a rod, in a vertical compartment. In this way, zero is represented by the flat position of the rod, while 9 is the highest position, at some distance from the top, which is a forbidden position (cipher positions are located at equally increasing levels). This is a ``zero altimetric abacus''. If otherwise, flat position is forbidden, but at same time, the top position is the highest possible one, we pass from the zero abacus to a zeroless abacus, and all the properties that we investigated in the paper easily translate in the different, but  related, ways of positioning rods in the two abaci. 

In the Middle Ages, positional (decimal) notation reached maturity and a decimal mark (a point or a comma) appeared in order to extend the decimal notion to fractions. In this way, following a different track, a denotation for zero naturally could arise. In fact, if the decimal mark is a point separating the integer from the proper fractional part, then a decimal point alone (without  preceding or following ciphers) is a denotation of zero. 

So far we discussed some historical implications of our analysis about positional systems (with and without zero). However, this analysis is not only relevant for the past mathematics. In fact, some aspects of number representation structures and algorithms are related to emergent aspects in discrete mathematics. Firstly, the algorithm outlined in Table \ref{symbolic} can be easily expressed as a particular {\it membrane system} (or P system) in the sense of \cite{paun}, where columns are membranes and the manipulation rules are transformation and movement of objects between membranes. In \cite{manca} (Chapter 2) an example is given of an abacus as a membrane system following the same idea, the multiplication method of Table \ref{symbolic} provides an improvement of the system in \cite{manca}, because computation can fully benefit from the distributed mechanism of membrane computing.  Moreover, the interest in string enumeration, as suggested by the examples of Section 1, is related to the mathematical analysis of genomes. They can be viewed as strings over the four letters $A,C,G,T$. When we fix a representation system of numbers, by means of sequences over the genomic alphabet, then genomes become numbers of gigantic sizes, therefore the methods of string (and number) representations have a direct impact for the investigations that aim at finding unconventional ways of considering genomes, and new perspectives to their study and understanding. 

\end{document}